\documentclass[12pt]{amsart}
\usepackage{amssymb}
\usepackage[all]{xy}

 \rm

\renewcommand{\(}{\left(}
\renewcommand{\)}{\right)}

\newcommand{\<}{\langle}
\renewcommand{\>}{\rangle}

\renewcommand{\bar}{\overline}
\newcommand{\abs}[1]{\left\lvert#1\right\rvert}
\newcommand{\norm}[1]{\left\lVert#1\right\rVert}
\newcommand{\st}{\:|\:}

\newcommand{\C}{{\mathbb{C}}}
\newcommand{\R}{{\mathbb{R}}}
\newcommand{\Z}{{\mathbb{Z}}}

\renewcommand{\phi}{\varphi}

\renewcommand{\H}{{\mathcal{H}}}

\theoremstyle{plain}
\newtheorem{thm}{Theorem}[section]
\newtheorem{lem}[thm]{Lemma}

\newtheorem{ex}[thm]{Example}
\newtheorem{conj}[thm]{Conjecture}

\theoremstyle{definition}

\theoremstyle{remark}

\title{Hermite expansions and Hardy's theorem}

\author{M.~K.~Vemuri}
\address{Chennai Mathematical Institute, Plot H1, SIPCOT IT Park, Padur~PO,
Siruseri 603103.}

\begin{document}

\begin{abstract}
Assuming that both a function and its Fourier transform are dominated
by a Gaussian of large variance, it is shown that the Hermite
coefficients of the function decay exponentially.  A sharp estimate
for the rate of exponential decay is obtained in terms of the
variance, and in the limiting case (when the variance becomes so small
that the Gaussian is its own Fourier transform), Hardy's theorem on
Fourier transform pairs is obtained.  A quantitative result on the
confinement of particle-like states of a quantum harmonic oscillator
is obtained.  A stronger form of the result is conjectured.  Further,
it is shown how Hardy's theorem may be derived from a weak version of
confinement without using complex analysis.
\end{abstract}

\maketitle

%\tableofcontents\newpage

\section{Introduction}\label{S:intro}

If $f \in L^1(\R)$, the {\em Fourier transform} of $f$ is defined by
$$
\hat f(\xi) = \frac{1}{\sqrt{2\pi}} \int f(x) e^{-i\xi x} \, dx.
$$
Let $g_a(x) = e^{-ax^2/2}$.  Hardy's theorem is usually stated as follows
(see \cite[Theorem 7.6]{Folland-Sitaram}, where the notation is slightly
different).

\begin{thm}\label{T:hardy}
For $a>0$, let
$$
E(a) = \{ f \in L^1(\R) \st
\abs{f(x)} \le C g_a(x) \text{ and } \abs{\hat f(\xi)} \le C g_a(\xi)
\quad\text{for some $C \in \R$} \}.
$$
If $a > 1$ then $E(a) = 0$.  If $a=1$ then $E(a) = \C g_a$.  If $a < 1$
then $\dim E(a) = \infty$.
\end{thm}

The last part of the trichotomy is usually substantiated by showing that
all Hermite functions belong to $E(a)$, if $a < 1$.

This statement of Hardy's theorem appears to suggest that if $a < 1$ then
no significant restriction is placed on $f$.  This is far from the truth.
In fact, regardless of the value of $a$, elements of $E(a)$ may be
characterized by the rate of exponential decay of their Hermite
coefficients.

%This paragraph is badly written.  Fix it later.
%
Hardy's theorem is usually proved by applying the Phragmen-Lindel\"of
principle to the Fourier-Laplace transform of $f$.  Instead, we apply
the Phragmen-Lindel\"of principle to the Bargmann transform (the
unitary intertwiner between the Schr\"odinger and Fock realizations of
the canonical representation of the Heisenberg group).  This transform
is better suited for studying Hermite expansions.

The result on exponential decay of Hermite coefficients leads, via
Mehler's formula to a Gaussian bound on the solutions of the Schr\"odinger
equation for the harmonic oscillator Hamiltonian, when the initial data
belong to $E(a)$.  We refer to this result as {\em confinement}.  We state
a stronger conjecture.

The solutions of the harmonic oscillator Schr\"odinger equation are
orbits of the standard maximal compact subgroup $K={\mathrm{SO}}(2)$
of ${\mathrm{SL}}(2, \R)$ under the metaplectic representation.
Further, the $K$-types are precisely the Hermite functions.  Using
this idea, we show that Hardy's theorem follows from a weak version of
confinement.  Thus, if a weak confinement result is proved by purely
PDE methods, we would have a proof of Hardy's theorem that does not
use complex analysis.  This would answer a question of Sundari.

Others have considered the connection between Hardy's theorem and
Schr\"odinger equations.  Chanillo \cite{Chanillo} showed that Hardy's
theorem is equivalent to a uniqueness theorem for the free-particle
Schr\"odinger equation.  The free-particle flow is the orbit of a
unipotent subgroup of ${\mathrm{SL}}(2, \R)$ under the metaplectic
representation.  It would be interesting to understand the connection
between Hardy's theorem and the metaplectic representation better;
perhaps there is a purely representation theoretic proof of Hardy's
theorem!

For more on the connections between analysis and the metaplectic
representation, see Howe \cite{Howe} or Folland \cite{Folland}.  For
more on Hardy's theorem, see Thangavelu \cite{Thangavelu}.
Information on Hermite functions and Mehler's formula may also be
found in \cite{Folland, Thangavelu}.

In this work, we will use the measure $dm = dx/\sqrt{2\pi}$ to define
the norm on $L^p(\R)$.

\section{Exponential decay of Hermite coefficients}\label{hermite}

We will use some properties of the Bargmann transform (see
\cite[p78]{Corwin-Greenleaf}, where there seems to be a normalization
error) in the proof of the main theorem.  To avoid cluttering up the
main argument, we recall these first.

Let $\H$ denote the Hilbert space of all entire functions $F$ on $\C$
such that
$$
\norm{F}^2 = \int \abs{F(w)}^2 \, \frac{e^{-\abs{w}^2/2} \, du \, dv}{2\pi}
 < \infty \quad\text{($w=u+iv$)}.
$$
Define $U:L^2(\R) \to \H$ by
$$
Uf(w) = \frac{e^{-w^2/4}}{2^{1/4}\pi^{1/2}} \int e^{xw} e^{-x^2/2} f(x) \, dx.
$$
Then $Uf$ is defined for Schwartz class functions $f$, and extends to an
isometric isomorphism.  We call $U$ the {\em Bargmann transform}.
Note that
$$
(U{\hat f})(w) = Uf(-iw),
$$
for all $w \in \C$.  Further, if $\phi_k$
denotes the $k$-th normalized Hermite function, then
$$
U\phi_k(w) = \frac{w^k}{\sqrt{2^k k!}}.
$$

\begin{thm}\label{T:exponential-decay}
Let $a \in (0, 1)$.  If
$$
\abs{f(x)} \le C g_a(x) \text{ and } \abs{\hat f(\xi)} \le C g_a(\xi)
$$
then
$$
\abs{\< f, \phi_k \>} \le
C \sqrt{\frac{2\pi k!}{1+a}} (e/k)^{k/2} \(\frac{1-a}{1+a}\)^{k/4}
$$
for $k=1,2,\dots$.
\end{thm}

\begin{proof}
Write $w=u+iv=re^{i\theta}$.  From the first hypothesis, we obtain
$$
\begin{aligned}
\abs{\int e^{xw} e^{-x^2/2} f(x) \, dx}
\le & \; C \int e^{xu - (1+a)x^2/2} \, dx \\
  = & \; C e^\frac{u^2}{2(1+a)} 
           \int e^{-\frac{1+a}{2} \(x - \frac{u}{1+a}\)^2} \, dx \\
  = & \; C \sqrt{\frac{2\pi}{1+a}} e^{\frac{u^2}{2(1+a)}}.
\end{aligned}
$$
Therefore,
$$
\begin{aligned}
\abs{Uf(w)}
\le & \; C \sqrt{\frac{2\pi}{1+a}}
           \exp\(\frac{v^2 - u^2}{4} + \frac{u^2}{2(1+a)}\) \\
  = & \; C \sqrt{\frac{2\pi}{1+a}}
           \exp{\frac{v^2 + \mu u^2}{4}} \\
  = & \; C \sqrt{\frac{2\pi}{1+a}}
           \exp \frac{(\mu + (1-\mu)\sin^2 \theta)r^2}{4},
\end{aligned}
$$
where $\mu=\frac{1-a}{1+a}$.

From the second hypothesis and the previous calculation, we obtain
$$
\begin{aligned}
\abs{Uf(w)}
  = & \; \abs{U{\hat f}(iw)} \\
\le & \; C \sqrt{\frac{2\pi}{1+a}}
           \exp \frac{(\mu + (1-\mu)\sin^2 (\theta+\pi/2))r^2}{4} \\
  = & \; C \sqrt{\frac{2\pi}{1+a}}
           \exp \frac{(\mu + (1-\mu)\cos^2 \theta)r^2}{4}.
\end{aligned}
$$

A substantial improvement in these estimates may be obtained by applying
the Phragmen-Lindel\"of principle to the {\em holomorphic} function $Uf$.
Let $\theta_0 = \frac{1}{2}\arctan\(\frac{2\sqrt{\mu}}{1-\mu}\)$,
$\theta_1=\frac{\pi}{2} - \theta_0$.  Observe that
$\theta_1 - \theta_0 < \frac{\pi}{2}$.  Let
$$
F(w) = \exp\(i \frac{\sqrt{\mu}}{4} w^2\) Uf(w).
$$
Then $F$ is entire, bounded by $3Ce^{\abs{w}^2}$ everywhere,
and by $C \sqrt{\frac{2\pi}{1+a}}$ on the rays $\theta=\theta_0$ and
$\theta=\theta_1$. It follows from the Phragmen-Lindel\"of principle that
$$
\abs{F(w)} \le C \sqrt{\frac{2\pi}{1+a}}
$$
for $\theta_0 \le \theta \le \theta_1$.  So
\begin{equation}\label{E:cask-strength}
\abs{Uf(w)} \le C \sqrt{\frac{2\pi}{1+a}}
\exp\(\frac{\sqrt{\mu}\sin 2\theta}{4} r^2\)
\end{equation}
for $\theta_0 \le \theta \le \theta_1$.  Combining this with the
previous two estimates, we obtain a crude estimate for $Uf$ in the
first quadrant:
$$
\abs{Uf(w)} \le C \sqrt{\frac{2\pi}{1+a}}
\exp\(\frac{\sqrt{\mu}}{4} r^2\).
$$
The same argument works in the other three quadrants, and so the
estimate holds everywhere.

If $Uf(w) = \sum_{n=1}^\infty c_n w^n$, the Cauchy estimates give
$$
\abs{c_n}
\le
C \sqrt{\frac{2\pi}{1+a}}
\exp\(\frac{\sqrt{\mu}}{4} r^2\) r^{-n}
$$
for all $r>0$.  Optimizing with respect to $r$, we get
$$
\abs{c_n}
\le
C \sqrt{\frac{2\pi}{1+a}}
\( \frac{e \sqrt{\mu}}{2n}\)^{n/2}.
$$
Therefore
$$
\begin{aligned}
\abs{\<f, \phi_k\>}
  = & \; \abs{\<Uf, U\phi_k\>} \\
  = & \; \int\int \(\sum_{n=0}^\infty c_n w^n\)
                  \bar{\(\frac{w^k}{\sqrt{2^k k!}}\)}
                  \frac{e^{-r^2/2} \, du \, dv}{2\pi} \\
  = & \; \frac{\abs{c_k}}{\sqrt{2^k k!}}
         \int\int r^{2k} \frac{e^{-r^2/2} \, du \, dv}{2\pi} \\
  = & \; \sqrt{2^k k!} \abs{c_k} \\
\le & \; C \sqrt{\frac{2\pi k!}{1+a}} (e/k)^{k/2} \mu^{k/4}
\end{aligned}
$$
\end{proof}

If $f \in E(1)$, then there exists a constant $C$ such that
$$
\abs{f(x)} \le C g_a(x) \text{ and } \abs{\hat f(\xi)} \le C g_a(\xi)
$$
for all $a \in (0,1)$.  So for $k \ge 1$ we have
$$
\abs{\<f, \phi_k\>}
\le C \sqrt{\frac{2\pi k!}{1+a}} (e/k)^{k/2} \mu^{k/4},
$$
for all $\mu \in (0,1)$.  It follows that $\<f, \phi_k\> = 0$ for
$k \ge 1$, and so $f \in \C \phi_0$.  If $a > 1$ and $f \in E(a)$,
then in particular $f \in E(1)$, and so $f = C \phi_0$.  However,
$\phi_0 \notin E(a)$, so $C=0$ and $f=0$.  So the classical Hardy
theorem follows from Theorem \ref{T:exponential-decay}.

If $a \in (0,1)$ and $f \in E(a)$ then
$$
\<f, \phi_k\> = O\(k^{1/4}\(\frac{1-a}{1+a}\)^{k/4}\)
$$
by Theorem \ref{T:exponential-decay} and the bound
$k! \le 3\sqrt{k}(k/e)^k$, $k=1,2,\dots$.  In particular, if
$a \in (0,1)$, $f \in E(a)$ and $\tanh(2\alpha) < a$ then
\begin{equation}\label{E:exponential-decay}
\<f, \phi_k\> = O(e^{-\alpha k}).
\end{equation}

To obtain the endpoint estimate ($\tanh(2\alpha) = a$), we need to use
the full strength of the estimate (\ref{E:cask-strength}).

\begin{thm}\label{T:sharp-exponential-decay}
If $f \in E(\tanh(2\alpha))$, then
$$
\<f, \phi_k\> = O(e^{-\alpha k}).
$$
\end{thm}

\begin{proof}
We will use the notation from the proof of Theorem \ref{T:exponential-decay}
with $a=\tanh 2\alpha$.  So $\mu=e^{-4\alpha}$.  Assume that
$f$ has norm at most $1$.  Define
$$
r_n(t)=
\begin{cases}
\sqrt{\frac{2n+2}{\mu + (1-\mu)\sin^2 t}}, & 0 \le t < \theta_0 \\
\sqrt{\frac{2n+2}{\sqrt{\mu} \sin 2t}}, & \theta_0 \le t \le \frac{\pi}{4}.
\end{cases}
$$
Extend $r_n$ to $[0,\pi/2]$ by the rule
$$
r_n(t) = r_n(\frac{\pi}{2} - t), \quad \frac{\pi}{4} < t \le \frac{\pi}{2},
$$
and to $[0,2\pi]$ by $(\pi/2)$-periodicity.  Then $r_n$ is positive,
continuous and piecewise smooth.  Put $\gamma_n(t) = r_n(t) e^{it}$.
Then each $\gamma_n$ winds once about the point $w=0$.
By the Cauchy integral formula, the estimate (\ref{E:cask-strength}) and the
eightfold symmetry,
$$
\begin{aligned}
\abs{c_n}
\le & \; \frac{1}{2\pi} \int_{\gamma_n} \abs{(Uf)(w)} \, \abs{w}^{-(n+1)} \,
                                        \abs{dw} \\
=   & \; \frac{4}{\pi} \sqrt{\frac{2\pi}{1+a}}
         \exp\(\frac{n+1}{2}\) (2n+2)^{-n/2} (I_n + J_n),
\end{aligned}
$$
where
$$
I_n = \int_0^{\theta_0} (\mu + (1-\mu)\sin^2 t)^{\frac{n-2}{2}}
                        \sqrt{\mu^2 + (1-\mu^2) \sin^2 t} \, dt
$$
and
$$
J_n = \mu^{n/4} \int_{\theta_0}^{\pi/4} (\sin 2t)^{\frac{n-2}{2}} \, dt.
$$
We estimate
$$
\begin{aligned}
I_n
\le & \; \int_0^{\theta_0} \(\frac{2\mu}{1+\mu}\)^{\frac{n-2}{2}}
                           \sqrt{\mu} \, dt \\
=   & \; \theta_0 \frac{1+\mu}{2\sqrt{\mu}} \(\frac{2\mu}{1+\mu}\)^{n/2},
\end{aligned}
$$
and
$$
\begin{aligned}
J_n
\le & \; \mu^{n/4} \int_0^{\pi/4} (\sin 2t)^{\frac{n-2}{2}} \, dt \\
=   & \; \frac{\sqrt{\pi}}{4}
         \frac{\Gamma\(\frac{n}{4}\)}{\Gamma\(\frac{n+2}{4}\)} \mu^{n/4}\\
\le & \; \frac{\sqrt{6\pi}}{4} n^{-1/2} \mu^{n/4}.
\end{aligned}
$$
Since $(2\mu)/(1+\mu) < \sqrt{\mu}$, it follows that $I_n = o(J_n)$,
and so
$$
c_n = O\(2^{-n/2} (e/n)^{n/2} n^{-1/2} \mu^{n/4}\).
$$
It follows, as before, that
$$
\<f, \phi_k\> = O\(k^{-1/4} \mu^{k/4}\)
              = O(e^{-\alpha k}).
$$
\end{proof}

\begin{ex}
With $a=\tanh(2\alpha)$, let
$$
f(x) = \exp\(\frac{-a + i \sqrt{1-a^2}}{2} x^2\).
$$
Then $f \in E(\tanh 2\alpha)$, but for all $\beta>1$ there exists
$c_\beta>0$ such that
$$
\abs{\<f, \phi_k\>} \ge c_\beta k^{-\beta/4} e^{-\alpha k},
\quad k=2, 4, 6, \dots.
$$
So Theorem \ref{T:sharp-exponential-decay} is sharp.
\end{ex}

\section{Confinement}\label{confinement}

The best constant $C$ in the definition of the space $E(a)$ (see Theorem
\ref{T:hardy}) is a norm on $E(a)$.  We won't introduce notation for it,
but will refer to it in context.  We find it convenient to reserve the
norm symbol for an $L^2$ type norm to be defined later.

Let $H=-\frac{\partial^2}{{\partial x}^2} + x^2$ denote the harmonic
oscillator.  Let $\psi_t(x)$ be a solution of the Schr\"odinger equation
\begin{equation}\label{E:hose}
\frac{1}{i} \frac{\partial \psi}{\partial t} = H \psi.
\end{equation}

\begin{thm}\label{T:confinement}
If $\psi_0 \in E(\tanh 2\beta)$ and $\gamma < \beta$ then for all $t \in \R$
$$
\psi_t \in E(\tanh \gamma),
$$
with bounded norm.
\end{thm}

The following proof was inspired by the proof of \cite[Theorem 9]{PSST}.

\begin{proof}
Assume $\psi_0 \in E(\tanh 2\beta)$ and $\gamma < \beta$.  Choose
$\gamma' \in (\gamma, \beta)$ and put $r=\gamma/\gamma'$.
Then $r \in (0,1)$.
The hypothesis and inequality (\ref{E:exponential-decay}) imply that
$$
\<\psi_0, \phi_k\> = O(e^{-\gamma' k}).
$$
If we write $\psi_0 = \sum_{n=0}^\infty \<\psi_0, \phi_n\> \phi_n$, then
$$
\psi_t = \sum_{n=0}^\infty e^{(2n+1)it}  \<\psi_0, \phi_n\> \phi_n.
$$
By the Cauchy-Schwarz inequality, and Mehler's formula
$$
\begin{aligned}
\abs{\psi_t(x)}
\le & \; \(\sum_{n=0}^\infty \abs{\<\psi_0, \phi_n\>}^{2(1-r)}\)^{1/2}
         \(\sum_{n=0}^\infty \abs{\<\psi_0, \phi_n\>}^{2r}
                             \abs{\phi_n(x)}^2\)^{1/2} \\
\le & \; \frac{1}{1-e^{-2(\gamma' - \gamma)}}
         \(\sum_{n=0}^\infty e^{-2\gamma n}
         \abs{\phi_n(x)}^2 \)^{1/2} \\
= & \; C(\gamma, \gamma') e^{-\frac{\tanh \gamma}{2} x^2}.
\end{aligned}
$$
Also,
$$
\abs{\widehat{\psi_t}(x)}
= \abs{\psi_{(t - \pi/4)}(x)}
\le  C(\gamma, \gamma') e^{-\frac{\tanh \gamma}{2} x^2}.
$$
So $\psi_t \in E(\tanh \gamma)$.
\end{proof}

We interpret Theorem \ref{T:confinement} as a result on the confinement
of particle-like states of the harmonic oscillator.  Regard the space
$E(a)$ (strictly speaking its image in projective space) as a
``Gaussian phase-box'' of side $1/a$.  If a state $\psi_0$ is initially
in the phase-box of side $\coth(2\beta)$ then its evolution $\psi_t$
is confined to the larger phase-box of side
$\coth(\beta-\varepsilon)$.

The following conjecture and example show that Theorem \ref{T:confinement}
is almost sharp.

\begin{conj}
If $\psi_0 \in E(\tanh 2\beta)$ then for all $t \in \R$
$$
\psi_t \in E(\tanh \beta).
$$
\end{conj}

The following example shows that we cannot do better.

\begin{ex}
Choose a branch $\sqrt{}$ of the square root that is defined on the right
half plane and is positive on the positive real line.  For
$\beta > 0$, let $r=e^{-2\beta}$, and
$$
\psi_{(t-\frac{\pi}{8})} =
\frac{e^{it}}{\sqrt{1+re^{4it}}}
\exp\(- \frac{1-re^{4it}}{1+re^{4it}} \frac{x^2}{2}\)
$$
Then $\psi$ is a solution of (\ref{E:hose}),
$$
\begin{aligned}
\abs{\psi_0}
=   & \; C_0 g_{\tanh(2\beta)} \\
\abs{\widehat{\psi_0}} = \abs{\psi_{-\frac{\pi}{4}}}
=   & \; C_{\frac{\pi}{4}} g_{\tanh(2\beta)}, \quad{\text{but}}\\
\abs{\psi_{-\frac{\pi}{8}}}
=   & \; C_{\frac{\pi}{8}} g_{\tanh(\beta)}.
\end{aligned}
$$
\end{ex}

\section{Confinement implies exponential decay}

In this section, we will show that the Hermite coefficients of a
``bound state'' decay exponentially.  We start with a simple estimate
for factorials that is slightly stronger than what can be obtained
from the standard Stirling formula.

\begin{lem}\label{L:factorial-estimate}
If $\beta > 1$ then there exists $B_\beta > 0$ such that
$$
2^{-2n}\frac{(2n)!}{(n!)^2} \ge B_\beta n^{-\beta/2}, \quad n=1,2,\dots.
$$
\end{lem}

\begin{proof}
Clearly, we need to prove this only for large $n$.  Note that
there exists $\delta>0$ such that $0 \le x \le \delta$ implies
$\log(1-x) \ge -\beta x$.  Choose $m$ so large that $k > m$
implies $0 \le \frac{1}{2k} \le \delta$.  Put
$$
D_\beta = \sum_{k=1}^m \log\(1-\frac{1}{2k}\),
\quad
B_\beta = e^{D_\beta} m^{\beta/2}.
$$
Let
$$
Q_n = 2^{-2n}\frac{(2n)!}{(n!)^2}.
$$
Then
$$
\begin{aligned}
\log Q_n
\ge & \; D_\beta - \frac{\beta}{2} \sum_{k=m+1}^n \frac{1}{k} \\
\ge & \; D_\beta - \frac{\beta}{2} (\log n - \log m).
\end{aligned}
$$
The result follows by exponentiation.
\end{proof}

  The results are most natural in an
$L^2$ setting.  So we define $E^2(a)$ to be the Hilbert space of all
functions $f$ such that
$$
2 \norm{f}_a^2 = \int \abs{f(x)}^2 e^{ax^2} \, \frac{dx}{\sqrt{2\pi}}
                 + \int \abs{\hat f(\xi)}^2 e^{a \xi^2}
                                          \, \frac{d\xi}{\sqrt{2\pi}}
             < \infty
$$
Observe that $a_1 < a_2$ implies $E(a_2) \subseteq E^2(a_1)$.

\begin{thm}\label{T:confinement-implies-exponential}
For all $\alpha>1/2$, there exists $A_\alpha > 0$ such that if
$a \in (0,1)$, and $\psi_t$ is a solution of (\ref{E:hose})
with $\norm{\psi_t}_a < C$ for all $t \in \R$ then
$$
\abs{\<\psi_0, \phi_k\>}
\le (C/A_\alpha) k^{\alpha/2} \(\frac{1-a}{1+a}\)^{k/2}
$$
for $k=1,2,\dots$.
\end{thm}

\begin{proof}
Let $f=\psi_0$, and for $n \in \Z$, let
$$
f_n = \frac{1}{2\pi}\int_0^{2\pi} \psi_t \bar{e^{int}} \, dt.
$$
Then
$$
\norm{f_n}_a \le C.
$$
%
%% By the Peter-Weyl theorem,
%% $$
%% f = \sum_{n \in \Z} f_n.
%% $$
%
Since $f_n$ is an eigenfunction of the harmonic oscillator with eigenvalue
$n$, we have $f_n=0$ if $n$ is even or negative, and
$$
f_{2k+1} = \<f, \phi_k\> \phi_k, \quad k=0,1,2,\dots.
$$
We will get a lower bound on $\norm{\phi_k}_a$.  This will imply an upper
bound on $\abs{\<f, \phi_k\>}$.

From Mehler's formula, we have
$$
\sum_{k=0}^\infty (\phi_k(x))^2 w^k
= \sqrt{2}(1-w^2)^{-1/2} e^{-\frac{1-w}{1+w} x^2}.
$$
Multiplying both sides by $e^{a x^2}$, integrating, and observing that
$\phi_k$ are real and are their own Fourier transforms, up to phase, we obtain
$$
\sum_{k=0}^\infty \norm{\phi_k}_a^2 w^k
 = (1-a)^{-1/2} (1-w)^{-1/2} (1 - w/\mu)^{-1/2},
$$
where $\mu=\frac{1-a}{1+a}$.  Expanding the right hand side in powers
of $w$, and equating coefficients, we obtain
$$
\norm{\phi_n}_a^2
 = (1-a)^{-1/2} 2^{-2n}
   \sum_{k=0}^n \frac{(2k)! (2(n-k))!}{(k! (n-k)!)^2} \mu^{-k}
$$
Since the above sum has non-negative terms, we must have
$$
\norm{\phi_n}_a^2 \ge (1-a)^{-1/2} 2^{-2n} \frac{(2n)!}{(n!)^2} \mu^{-n}.
$$
So by Lemma \ref{L:factorial-estimate} if $\alpha > 1/2$ there exists
a constant $A_\alpha>0$ such that
$$
\norm{\phi_k}_a \ge A_\alpha (1-a)^{-1/4} k^{-\alpha/2} \mu^{-k/2},
\quad k=1,2,\dots,
$$
and so
$$
\abs{\<f, \phi_k\>} \le (C/A_\alpha) (1-a)^{1/4} k^{\alpha/2} \mu^{k/2}
                      = (C/A_\alpha) (1-a)^{1/4} k^{\alpha/2}
                        \(\frac{1-a}{1+a}\)^{k/2}, \quad k=1,2,\dots.
$$

\end{proof}

%% We can derive the Cowling-Price $L^2$-version of Hardy's theorem from
%% Theorem \ref{T:confinement-implies-exponential}.
%% If $a=1$ this implies $f=0$.

Theorem \ref{T:confinement-implies-exponential} suggests a new approach
to proving Hardy's theorem.  Using PDE methods, we first prove

\begin{thm}[Weak confinement]\label{T:weak-confinement}
There exist $N$ such that for all $\beta>0$, if $\psi_t$ is a solution
of (\ref{E:hose}) and $\psi_0 \in E^2(\tanh(N\beta))$ then there exists
$K$ such that
$$
\norm{\psi_t}_{\tanh \beta} \le K \norm{\psi_0}_{\tanh(N\beta)}
$$
for all $t \in \R$.
\end{thm}

Write $a=\tanh \beta$ and $b=\tanh(N\beta)$.  If $\psi_0 \in E(1)$ with
norm bounded by $1$, then
$$
\norm{\psi_0}_b \le 2^{-1/4}(1-b)^{-1/4}
$$
for all $\beta>0$.  So by Theorem \ref{T:weak-confinement}, there exists
$K$ such that
$$
\norm{\psi_t}_a \le K (1-b)^{-1/4}
$$
So by Theorem \ref{T:confinement-implies-exponential}, there exists
$A>0$ such that for all $\beta>0$ we have
$$
\begin{aligned}
\abs{\<\psi_0, \phi_k\>}
\le & \; \frac{K(1-b)^{-1/4}}{A} (1-a)^{1/4} k \(\frac{1-a}{1+a}\)^{k/2} \\
\le & \; \frac{Kk}{A} e^{\frac{(N-1)\beta}{2}} e^{-\beta k} \\
=   & \; \frac{Kk}{A} e^{\beta \(\frac{N-1}{2} - k\)}
\quad k=1,2,\dots
\end{aligned}
$$
It follows that $\psi_0$ is a {\em finite}
linear combination of Hermite functions.  Since
$\psi_0 \in E(1)$, it follows that the corresponding linear
combination of Hermite polynomials is bounded, and hence constant.
So $\psi_0$ is a Gaussian.

\bibliographystyle{amsplain}
\bibliography{v10-heht}

\end{document}